\documentclass[12pt,twoside]{article}
\usepackage{amssymb}
\usepackage{amsmath,amscd}
\usepackage{pst-node,pstricks,multido,pst-plot,pst-text,pst-3d}%
\usepackage{graphicx}
\usepackage{amsfonts}

\def\shuf{\sqcup\!\!\sqcup}
\def\square#1{\hbox{\vrule width \thickness
   \vbox to \squaresize{\hrule height \thickness\vss
      \hbox to \squaresize{\hss#1\hss}
   \vss\hrule height\thickness}
\unskip\vrule width \thickness} \kern-\thickness}

\def\vsquare#1{\vbox{\square{$#1$}}\kern-\thickness}

\def\gfill{\leaders\hbox to 1.2em{\hss\g\hss}\hfill}

\def\tab#1|#2|{ {}_{#2}^{#1} }   

\def\ttab#1|#2|#3|#4|{
\vbox{\hbox{$\tab#1|#2|$} \hbox{$\tab#3|#4|$} } }

\font\goth=eufm10 scaled 1200
\def\hfl#1#2{\smash{\mathop{\hbox to 12mm{\rightarrowfill}}
\limits^{\scriptstyle#1}_{\scriptstyle#2}}}

\def\build#1_#2^#3{\mathrel{
\mathop{\kern 0pt#1}\limits_{#2}^{#3}}}

 \font\goth=eufm10 scaled \magstep1

\newtheorem{example}{Example}[section]

\newtheorem{proposition}[example]{Proposition}

\def\adots{\mathinner{\mkern2mu\raise1pt\hbox{.}
\mkern3mu\raise4pt\hbox{.}\mkern1mu\raise7pt\hbox{.}}}
\def\<{\langle}
\def\>{\rangle}

\def\cf{{\it cf.}}
\def\shuf{\mathrel{\sqcup}\joinrel\mathrel{\sqcup}}

\def\ie{{\it i.e. }}

\def\Sym{{\bf Sym}}

\def\N{\mathbb N}

\def\F{\mathbb F}
\def\M{\mathbb M}

\def\Q{\mathbb Q}

\newtheorem{theo}{Theorem}
\newtheorem{prop}[theo]{Proposition}

\newtheorem{cor}[theo]{Corollary}

\newtheorem{ex}{Example}

\begin{document}

\title
{NONCOMMUTATIVE SYMMETRIC FUNCTIONS ASSOCIATED WITH A CODE, LAZARD
ELIMINATION, AND WITT VECTORS }
\author{ \sc Jean-Gabriel \sc Luque \ and  Jean-Yves \sc Thibon\\ \\
IGM,   Laboratoire d'informatique\\ 77454 Marne-la-Vall\'ee Cedex 2,
France.\\ \{luque,jyt\}@univ-mlv.fr}
 \maketitle
  \begin{abstract} The
construction of the universal ring of Witt vectors is related to
Lazard's factorizations of free monoids by means of a noncommutative
analogue. This is done by associating to a code a specialization of
noncommutative symmetric functions.
\end{abstract} {\bf Keywords:} {Witt vectors, codes, symmetric
functions, factorizations of free monoids}

\section{Introduction}

Among Ernst Witt's many contributions to mathematics, one finds
two apparently unrelated ideas, both published in 1937 in two
consecutive issues of Crelle's journal. The first one, the ring of
Witt vectors \cite{Witt1} (a generalisation of $p$-adic numbers),
solves a problem in commutative algebra, whilst the second one,
the introduction of the free Lie algebra \cite{Witt2}, definitely
pertains to  noncommutative mathematics.

The aim of this note is to point out a close connection between
both constructions, through the notion of noncommutative symmetric
functions. What comes out is that the natural noncommutative
analogues of the symmetric functions classically associated to the
construction of universal Witt vectors can be related to another
classical topic in the combinatorics of free Lie algebras:
Lazard's factorizations of free monoids. This relation manifests
itselfs when the elementary noncommutative symmetric functions are
specialized by means of a code.

Our notations will be those of \cite{Mcd} and \cite{GKLLRT}. If
$S$ is a set of words,  we denote by $\underline S=\sum_{w\in S}w$
its characteristic series.

\section{Noncommutative Witt symmetric functions}

Let us first recall the construction of the universal ring $W(R)$
of Witt vectors associated to a commutative ring $R$.  This is not
Witt's original construction, but a simpler one that he found much
later, and that he communicated to S. Lang, who published it as a
series of exercises in his famous textbook of Algebra \cite{Lang}.

The ring $W(R)$ can be characterized by the following properties
\cite{Car} (see also \cite{Ha,Di}):

\bigskip (W1)  As a set, $W(R)=\{{\bf a}=(a_n)_{n\ge 1}\, ,\, a_n\in R\}$,
and for any ring homomorphism $f:\ R\rightarrow R'$, the map
$W(f):\ {\bf a}\mapsto (f(a_n))_{n\ge 1}$ is a ring homomorphism.

\medskip (W2) The maps $w_n:\ {\bf a}\mapsto \sum_{d|n}da_d^{n/d}$ are
ring homomorphisms $W(R)\rightarrow R$.

\bigskip
Operations on Witt vectors are better understood in terms of
symmetric functions. Let $X=\{x_n\, ,\, n\ge 1\}$ be an infinite
set of commuting indeterminates (called here an {\it alphabet}),
and following \cite{Re1}, define symmetric functions $q_n(X)$ by
\begin{equation}\label{DEFq}
\prod_{n\ge 1} {1\over 1-t^n q_n(X)} = \sigma_t(X) :=\sum_{n\ge 1}
t^n h_n(X)
\end{equation}
where the complete homogeneous functions $h_n(X)$ are defined by
$\sigma_t(X)=\prod_{n\ge 1}(1-tx_n)^{-1}$ ({\it cf.} \cite{Mcd}).
The $q_n$'s are connected to the power-sums $p_n(X)=\sum_i x_i^n$
by
\begin{equation}\label{pc}
p_n = \sum_{d|n} dq_d^{n/d}
\end{equation}
and condition (W2) can be regarded as expressing the familiar
properties of power sums $$ p_n(X+Y)=p_n(X)+p_n(Y) $$ $$
p_n(XY)=p_n(X)p_n(Y) $$ where we use the $\lambda$-ring notation
(an alphabet is identified with the formal sum of its elements).
The transformation (\ref{DEFq}) is also used in \cite{BBP,MR}.

Reutenauer \cite{Re1} studied the symmetric functions $q_n(X)$ and
made the conjecture that for $n\ge 2$, $(-q_n)$ is Schur positive.
This conjecture was proved in \cite{ST2}, in a generalised
noncommutative form. Denoting as in \cite{GKLLRT} the
noncommutative complete homogeneous symmetric functions of an
alphabet $A$ by $S_n(A)$, one introduces noncommutative Witt
symmetric functions $Q_n(A)$  by the identity
\begin{equation}\label{DEFQ}
\prod_{n\ge 1}^{\rightarrow} {1\over 1-t^n Q_n(A)} = \sigma_t(A)
:=\sum_{n\ge 1} t^n S_n(A)
\end{equation}
Then, it is proved in \cite{ST2} that for $n\ge 2$, $(-Q_n)$ is a
positive, multiplicity free,  sum of noncommutative ribbon Schur
functions $R_I$.
\section{Noncommutative symmetric functions associated to a
code}\label{code} Let $A$ be an alphabet and ${C}\subset {A}^+$ a
{\it code}, \ie a minimal generating set of  a free submonoid
${C}^*$ in ${A}^*$. Suppose we have a decomposition of $C$ as a
disjoint union
\begin{equation}
{C}=\coprod_{n\ge 1} {C}_n
\end{equation}
of possibly empty subsets ${C}_n$. In what follows, we shall
in general assume that ${C}_n={C}\cap {A}^n$, but this restriction in not
necessary. In any case, we will consider that the elements of
${C}_n$ have degree $n$.

We will denote by $l(w)$ the length of a word $w$. Recall that, as
a graded algebra, ${\bf
Sym}=\Q\langle\Lambda_1,\Lambda_2,\dots\rangle$ where $\Lambda_n$
has degree $n$. We can  define a specialization $\Sym[C]$ of
the algebra of noncommutative symmetric functions by setting
\begin{equation}
\Lambda_n[C]=(-1)^{n-1}\underline{C}_n\,.
\end{equation}
With this choice of signs, the complete symmetric functions are
then given by
\begin{equation}
S_n[C]=\sum_{w\in (C^*)_n}w\, ,
\end{equation}
the sum of all elements of degree $n$ in $C^*$. We call  $C$-{\it
Witt symmetric function} the value $Q_i[C]$ of
 $Q_i$ under this specialization.
Let us give the first $Q_i[C]$ for three examples

\begin{ex}{\rm The Fibonacci prefix code ${C}=\{b,ab\}$.
}
\begin{equation}
\begin{array}{l}
Q_1[C]=b\\ Q_2[C]=ab\\ Q_3[C]=ab^2\\ Q_4[C]=ab^3\\
Q_5[C]=ab^2ab+ab^4\\ Q_6[C]=ab^3ab+ab^5\\
Q_7[C]=ab^6+ab^2(ab)^2+ab^4ab+ab^3ab^2
\end{array}\nonumber\end{equation}
\end{ex}

\begin{ex}\label{ex2}{\rm The infinite prefix code ${C}=ba^*$.
}
\begin{equation}
\begin{array}{l}
Q_1[C]=b\\ Q_2[C]=ba\\ Q_3[C]=ba^2+bab\\ Q_4[C]=bab^2+ba^3+ba^2b\\
Q_5[C]=bab^2a+bab^3+ba^2b^2+ba^2b^2+ba^4+ba^2ba\\
Q_6[C]=ba^4b+ba^3b^2+ba^2b^2a+bab^3a+bab^4+ba^3ba+ba^5+ba^2b^3
\end{array}\nonumber\end{equation}
\end{ex}

\begin{ex}\label{Dyck}{\rm The Dyck code $ D$ (for the Dyck language
$\underline{D^*}=1+a\underline{D^*}b\underline{D^*}.$) }
\begin{equation}
\begin{array}{l}
Q_1[D]=Q_3[D]=Q_5[D]=Q_7[D]=0\\ Q_2[D]=ab\\ Q_4[D]=aabb\\
Q_6[D]=aaabbb+aababb+aabbab
\\
Q_8[D]=aaaabbbb+aaababbb+aabaabbb+aaabbbab+aaabbabb\\
+aabbabab+aabababb+aababbab
\end{array}\nonumber\end{equation}
\end{ex}
On these examples, we remark that each $Q_i$ is  multiplicity free
and is the characteristic series of a code. In the following
section, we prove that it is always the case and  give a
characterization of these codes in terms of Lazard's
factorizations.
\section{Witt symmetric functions and factorization}
\subsection{Lazard elimination process}
We recall that a factorization of a monoid $\M$ is an ordered
family of monoids $\F=(\M_i)_{i\in I}$ such that each element
$m\in\M$ admits an unique decomposition
\begin{equation}
m=m_{i_1}\cdots m_{i_k}
\end{equation}
where $i_1>\cdots>i_k$ and
$m_{i_1}\in\M_{i_1},\dots,m_{i_k}\in\M_{i_k}$. In the case of free
monoids ${\M}=A^*$, this property can be stated in terms of
generating series
\begin{equation}\label{fact1}
\sum_{w\in {A}^*}w=\prod^{\leftarrow}_i\sum_{w\in \M_i}w
\end{equation}
A submonoid $\M'\subset \M$ can be
characterized by  its generating set $\M'\setminus\M'^2$. In
the sequel, a factorization will be denoted by the sequence of the
generating sets of its components, and for a factorization
${\F}=({C}_i)_{i\in I}$ of ${A}^*$ into free submonoids,
(\ref{fact1}) reads
\begin{equation}
\frac1{1-\underline A}=\prod^{\leftarrow}_i\frac1{1-\underline
C_i}.
\end{equation}
Factorisations of  free monoids have been extensively studied and
the reader can refer to \cite{BP,Vi} for a survey. A simple but
relevant example is  a {\it Lazard bisection}: considering a
subalphabet $B\subset A$, one has
\begin{equation}\label{bisection}
{A}^*={B}^*(({A}\setminus {B}){B}^*)^*.
\end{equation}
The pair $({B},({A}\setminus {B}){B}^*)$ is
a factorization of ${A}^*$ (\cf \cite{Vi,BP}, and \cite{Reut} for
applicationa to free Lie algebras).
 Now, from
(\ref{bisection}) we can obtain a trisection (\ie a factorization
in three submonoids) by iterating the process on the left or the
right factor, and so on. Factorisations which can be generated
applying only Lazard bisections on the right factors will be
called {\it Lazard right compositions}. This is a simple
particular  case of a {\it locally finite right regular
factorization} \cite{Vi}. Let $A$ be an alphabet (possibly
 infinite) and
$\omega: {A}\rightarrow \N\setminus0$ a weight function. This
function can be extended uniquely  as a morphism $\omega:
{A}^*\rightarrow (\N,+)$. We can associate a factorization
$\F({A},\omega)$ to this weight as follows. Let $({Z}_i)_{i\ge 1}$
and $({C}_i)_{i\ge 1}$  be the sequences  of codes defined by the
recurrence relations
\begin{enumerate}
\item ${Z}_1={A}$
\item $\mbox{for each integer }i>0, {C}_i={Z}_i\cap\{w\in
A^*|\omega(w)=i\}.$
\item $\mbox{for each integer }i>0, {Z}_{i+1}=({Z}_i\setminus {C}_i){C}_i^*$
\end{enumerate}
The sequence $\F({A},\omega)=({C}_i)_{i\in I}$ obtained by
omitting the empty elements in $({C}_i)_{i\ge 1}$ is a Lazard
right composition.  Since each code ${C}\in {A}^*$ admits the
length of words as  weight function, we define the {\it right
length factorization} $\F({C^*})$ of  $C^*$ by
\begin{equation}
\F({C^*})=\F({C},l).
\end{equation}
We can remark that a code is homogeneous if and only if
$\F({C^*})=({C})$.

\subsection{Computation of the $C$-Witt symmetric functions}
The equality between formal series (\ref{DEFQ}) can be rewritten
as
\begin{equation}
\frac1{1-\underline C} = \frac1{1-Q_1}\frac1{1- Q_2}\frac1{1-
Q_3}\cdots,
\end{equation}
and in this section,  we  prove that this
factorization of series is the weight right factorization of
${C}^*$.

In a more general setting, T. Scharf and one of the
authors \cite{ST2}, gave a recursive algorithm for computing the
$Q_n's$. Recall that the algebra {\bf Sym} of  noncommutative
symmetric functions is the free associative
algebra $\Q\langle S_1,S_2,\cdots \rangle$ generated by an
infinite sequence of non commuting variables $S_n$, graded by the
weight function $w(S_n)=n$. If $I=(i_1,\cdots,i_n)$, one defines
\begin{equation}\label{stilde}
\tilde S^I=(-1)^{n}S_{i_1}\cdots S_{i_n}
\end{equation}
The $Q_i$ can be computed following the rules given in \cite{ST2}:
\begin{enumerate}
\item $F_1=-\sum_{i} \tilde S_i$ ,
\item $F_{n+1}=F_n+Q_n(1-F_n)$,
\item $Q_n$ is the term of weight $n$ in $F_n$ multiplied by $-1$.
\end{enumerate}
Setting $Z_n=1-(1-F_n)^{-1}$, we obtain
\begin{equation}
F_{n+1}=1-(1-Z_{n+1})^{-1}
\end{equation}
and \begin{equation}
\begin{array}{rcl}
F_n+Q_n(1-F_n)&=&1-(1-Z_n)^{-1}+Q_n(1-Z_n)^{-1}\\
&=&1-(1-Q_n)(1-Z_n)^{-1}.
\end{array}
\end{equation}
This implies
\begin{equation}
\begin{array}{rcl}
Z_{n+1}&=&1-(1-Z_n)(1-Q_n)^{-1}\\ &=&(Z_n-Q_n)(1-Q_n)^{-1}.
\end{array}
\end{equation}
Following \cite{ST2}, each $Q_i$ is multiplicity free on the
$\tilde S_I$.

Let $Z_i[C]$ and  $F_i[C]$ the values of $Z_i$ and $F_i$ under the
specialization $S_n=S_n[C]$. The definition of noncommutative
complete functions gives
\begin{equation}
\sigma_1[C]=\sum S_n[C]=\frac1{1-\underline C}.
\end{equation}
Hence, one has
\begin{eqnarray}
Z_1[C]=1-(\sigma_1[C])^{-1}=\underline C,\\
Q_1[C]=S_1[C]=\underline{C_1},\mbox{ and }\\
Z_1[C]-Q_1[C]=\underline C-\underline{C_1}.
\end{eqnarray}
It follows that $Z_1[C]$ and $Q_1[C]$ are the characteristic
series of the  codes ${Z}_1={C}$ and ${Q_1}={C_1}$. By induction,
for each $n>0$ the series $Z_{n+1}[C]$ and $Q_{n+1}[C]$ are the
characteristic series of ${Z}_{n+1}=({Z_n}\setminus {Q}_n){Q}_n^*$
and ${Q}_n^*={Z_n}\cap A^{\leq n}$. Hence, the following statement
holds.
\begin{prop}
Let $C$ be a code. Each $C$-Witt symmetric function $Q_i[{C}]$
 is  the characteristic series of a code,  and the
sequence  obtained by deleting the empty set
from $({Q}_1,{Q}_2,\dots,{Q}_n,\dots)$ is the  right length
factorization of ${C}^*$.
\end{prop}
\begin{ex}\label{ex4}
{\rm The sequence $(a,Q_1[ba^*],Q_2[ba^*],\dots)$ is a
factorization of $A^*$. The same method is applicable to compute
homogeneous factorizations of non-homogeneous alphabets. For
example, considering the alphabet $A=\N\setminus 0$ with the
weight $\omega=id$, one finds
\begin{equation}\begin{array}{l}
Q_1[A]=1\\ Q_2[A]=2\\ Q_3[A]=21+3\\ Q_4[A]=211+31+4\\
Q_5[A]=2111+212+311+32+41+5\\
Q_6[A]=21111+51+2112+6+3111+312+42+411
\end{array}\nonumber\end{equation}
It is easy to see that this can be obtained from example \ref{ex2}
by the morphism $ba^n\rightarrow n+1$.}
\end{ex}
One has the following decomposition
\begin{cor}
Let $w\in  C^*$, then either $w\in C$ either there exist $w_1,
w_2\in  C^*$ such that $w=w_1w_2$, with $\omega(w_2)< \omega(w_1) <
\omega(w)$. Furthermore, if $w_1\not\in \ C$ then $w_1=w'_1w''_1$
with $w'_1, w''_1\in  C^*$ and $\omega(w''_1)\leq \omega(w_2)$
\end{cor}
This follows from the standard bracketing process of
regular factorizations, which is described in \cite{Vi}.

\subsection{Noncommutative elementary symmetric functions and
Lazard elimination}

The link between Lazard elimination and
noncommutative Witt symmetric functions can be better understood
in terms of elementary symmetric functions. The generating
function of noncommutative elementary symmetric function is
\begin{equation}
\lambda_t=\sum_{k\ge 0}\Lambda_kt^k.
\end{equation}
These  are related to noncommutative complete functions
by
\begin{equation}
\sigma_t=\frac1{\lambda_{-t}}=\frac1{1-\Lambda_1t+\Lambda_2t^2-\cdots+(-1)^n\Lambda_nt^n+\cdots}.
\end{equation}
If we set $\tilde \Lambda_n=(-1)^{n+1}\Lambda_n$, the series
$\sigma_1$ can be considered as the characteristic series of
the free monoid ${\bf \Lambda}^*=\{\tilde\Lambda_1,
\tilde\Lambda_2,\dots,\tilde \Lambda_n\,\dots\}^*$. We endow this monoid
with the weight function defined by $\omega(\tilde\Lambda_n)=n$.
Then,
\begin{theo}
One has
\begin{equation}
\F({\bf \Lambda},\omega)=(Q_1[{\bf \Lambda}],Q_2[{\bf
\Lambda}],\dots)=(Q_1,Q_2,\dots)
\end{equation}
\end{theo}
This provides a simple algorithm for computing the decomposition of
the  noncommutative symmetric Witt functions on the basis of
elementary symmetric functions. Let us gives the computation of
the first $Q_i$'s,
\begin{equation}\nonumber\begin{array}{ll}
Q_1=&\Lambda_1\\ Q_2=&-\Lambda_2\\
Q_3=&-\Lambda_2\Lambda_1+\Lambda_3\\
Q_4=&-\Lambda_2\Lambda_1\Lambda_1+\Lambda_3\Lambda_1-\Lambda_4\\
Q_5=&-\Lambda_2\Lambda_1^3+\Lambda_2\Lambda_1\Lambda_2+
\Lambda_3\Lambda_1\Lambda_1-\Lambda_3\Lambda_2-\Lambda_4\Lambda_1+\Lambda_5\\
Q_6=&-\Lambda_2\Lambda_1^4+\Lambda_5\Lambda_1+
\Lambda_2\Lambda_1\Lambda_1\Lambda_2-\Lambda_6+\Lambda_3\Lambda_1^3-\Lambda_3\Lambda_1\Lambda_2+
\Lambda_4\Lambda_2\\&-\Lambda_4\Lambda_1\Lambda_1\\
Q_7=&-\Lambda_3\Lambda_1^2\Lambda_2+
\Lambda_3\Lambda_1^4-\Lambda_3\Lambda_1\Lambda_2\Lambda_1+
\Lambda_3\Lambda_1\Lambda_3+\Lambda_3\Lambda_2^2-\Lambda_4\Lambda_1^3
\\&+\Lambda_4\Lambda_1\Lambda_2
-\Lambda_4\Lambda_3+\Lambda_4\Lambda_2\Lambda_1+\Lambda_5\Lambda_1^2-\Lambda_5\Lambda_2
-\Lambda_6\Lambda_1+\Lambda_7\\&-\Lambda_2\Lambda_1^5
+\Lambda_2\Lambda_1^3\Lambda_2+
\Lambda_2\Lambda_1^2\Lambda_2\Lambda_1-\Lambda_2\Lambda_1^2
\Lambda_3-\Lambda_2\Lambda_1\Lambda_2^2
\end{array}
\end{equation}
(compare  Examples \ref{ex2} and \ref{ex4}).

The decomposition of the elementary functions on the basis
$Q_I=Q_{i_1}\cdots Q_{i_n}$ is obtained by inspection of the
series
\begin{equation}
\frac1{\sigma_t}=\lambda_{-t}=1-\sum_n\tilde\Lambda_nt^n=\prod^{\leftarrow}_i(1-Q_it^i).
\end{equation}
One finds
\begin{equation}
\tilde\Lambda_n=\sum_{k}(-1)^k\sum_{i_1>\cdots>i_k\atop
i_1+\cdots+i_k=n}Q_{i_1}\cdots Q_{i_k}
\end{equation}
\section{Examples involving Gaudier's $*$-multinomials}
\subsection{Gaudier's $*$ multinomials}

Classical (commutative) Witt vectors give rise to
examples involving interesting integer sequences.
Let us recall for example the construction of $*$-multinomial
coefficients and $*$-factorials
given by Gaudier in \cite{Gaud}.
If $R$ is a $\Q$-algebra, there is a commutative diagram
\begin{eqnarray}\label{diag}
\nonumber\begin{CD} W({ R}) @>e>> \Lambda({ R})=1+t{ R}[[t]]\\ @V{w}VV
       @VV{\partial}V\\
{ R}^{\N^*} @>\iota>> { R}[[t]]
\end{CD}
\end{eqnarray}
where
\[
\partial=\frac d{dt}\ln,
\]
\[l(c_1,\cdots,c_n,\cdots)=
\sum_{n\geq 1}c_nt^{n-1},
\]
\[
w=(w_1,\cdots,w_n,\cdots),\]
\[
e(a_1,\dots,a_n,\dots)=\displaystyle\prod^\rightarrow\frac1{1-a_nt^n}.
\]
These maps are all isomorphisms. Let $i_1, i_2, \dots, i_k$ be $k$
positive integers and $n=i_1+\cdots+i_k$. In \cite{Gaud} Gaudier
has defined the $*$-multinomial coefficient  $*\left(n\atop
i_1,\dots,i_k\right)\in W(\Q)$ as the Witt vector such that
$w_p=\left(np\atop pi_1,\dots,pi_k\right)$. In the same paper, he
has defined the $*$-factorial $*n!/n!$ by
$w_p(*n!/n!)=\frac1{n!}\frac{(pn)!}{p!^n}$. In particular, he has
computed $e(*2!/2!)=e(\frac12*\left(2\atop1\right))$ in closed
form and related it to the Catalan numbers
\begin{equation}\label{f2!}
e(*2!/2!)=\left(1+\sqrt{1-4t}\over2\right)^{-1}
\end{equation}
This raised the question whether it was possible to find
similar expressions for
other $*$-factorials and
$*$-multinomials, and to give combinatorial
interpretations.

\subsection{Non-commutative analogues of
$\frac1n*\left(n\atop1\right)$\label{(n,1)/n}} Consider the Dyck
code of example \ref{Dyck}. The free monoid $D^*$ is a submonoid
of $\{aa,ab,bb\}^*$, graded by $\rho(aa)=\rho(bb)=\rho(ab)=1$.
Under the $D$-specialization,
\begin{equation}
\Lambda_n=(-1)^{n+1}\sum_{w\in D\atop\rho(w)=n}w.
\end{equation}
Recall that the noncommutative power sums of the second kind
$\Phi_n$ are defined by
\begin{equation}
\sum_{n\ge 1}\Phi_n\frac{t^n}{n}=\log \sigma_t\,.
\end{equation}
If one interprets Dyck words as binary trees, the $\Phi_n$ can be
regarded as sums over forests
\begin{equation}
\Phi_n=n\sum_k\frac1k\sum_{w_1,\dots,w_k\in D\atop
\rho(w_1)+\cdots+\rho(w_k)=n}w_1\cdots w_k.
\end{equation}
\begin{example}\rm  The first values of $\Phi_n$ are

$\Phi_1=ab$

$\Phi_2=2aabb+ab.ab$

$\Phi_3=3aaabbb+3aababb+\frac32aabb.ab+\frac32ab.aabb+ab.ab.ab$

$\begin{array}{rcl}\Phi_4&=&4aaaabbbb+4aaabbabb+4aabaabbb+4aabababb+4aaababbb\\
&&+2aaabbb.ab+2aababb.ab+2ab.aaabbb+2ab.aababb+\frac43aabb.ab.ab\\
&&+\frac43 ab.aabb.ab+\frac43 ab.ab.aabb+ab.ab.ab.ab
\end{array}$\\
See example \ref{Dyck} for the first values of $Q_n$.
\end{example}
Setting $\pi(w)=1$ for all words $w$, one has
\begin{equation}
\pi(\Phi_n)=w_n(\frac1{2!}*2!)=w_n(\frac1{2}*\left(2\atop
1\right))
\end{equation}
and one recovers (\ref{f2!}).

More generally, for
$*\frac1n\left(n\atop1\right)$, an easy application of the
Lagrange inversion formula gives that
\begin{equation}\label{n!}
e\left(\frac1n*\left(n\atop1\right)\right)=\sum_{k\ge
0}{\left(nk\atop k\right)\over (n-1)k+1}t^k
\end{equation}
is also the generating series of $n$-ary trees. Consider the free
monoid ${\goth F}_n=T_n^*$ over the alphabet $T_n=\{a_{\bf t}\}$
whose letters are labelled by $n$-ary trees and weighted by
$\rho(a_{\bf t})=E({\bf t})/n$, where
$E({\bf t})$ is the  number of edges of $\bf t$.
Instead of specializing $\Lambda_k$ as in the case  $n=2$,
we set $S_k=\sum_{\rho({\bf t})=k}a_{\bf t}$. Then, we  obtain
\begin{equation}
\Phi_k=k\sum_p{(-1)^{p+1}\over p}\sum_{{\bf t}_1,\dots, {\bf
t}_p\atop\rho({\bf t}_1)+\cdots+\rho({\bf t}_p)=k}a_{{\bf
t}_1}\cdots a_{{\bf t}_p}.
\end{equation}
Again, applying the morphism $\pi(w)=1$, we obtain
\begin{equation}
\pi(\Phi_k)=w_k(\frac1n*\left(n\atop1\right)),
\end{equation}
that is, we recover (\ref{n!}). Remark that this implies the identity
\begin{equation}
\sum_p{(-1)^{p+1}\over
p}\sum_{i_1+\cdots+i_p=k}\prod_m{\left(ni_m\atop i_m\right)\over
(n-1)i_m+1}=\frac1{nk}\left(nk\over k\right).
\end{equation}

\subsection{Combinatorial interpretations of some $*$-binomial coefficients}

Remark that, as in section \ref{(n,1)/n},
the specialization $S_k=\sum_{\rho({\bf t})=k+1}a_{\bf t}$ gives a
non-commutative analogue of
\begin{equation}
e(*\left(n\atop1\right))(z)=\sum_{k=1}^\infty{\left(nk\atop
k\right)\over
(n-1)k+1}t^{k-1}={e(\frac1n*\left(n\atop1\right))-1\over z}.
\end{equation}
This last equality can be generalized as follows.

\begin{proposition} Let $p\ge 1$ and $\omega=e^{2i\pi/p}$.
\begin{eqnarray}
e(*\left(np\atop p\right))&=&
\frac{-1}z\prod_{k=0}^{p-1}\left(1-e(\frac1n*\left(n\atop1\right)\right)
(\omega^k z^{1\over p})) \label{e1}\\
&=&\prod_{k=0}^{p-1}e\left(*\left(n\atop1\right)\right)
(\omega^kz^{1\over p})\label{e2}
\end{eqnarray}
\end{proposition}
{\bf Proof} We first prove (\ref{e2}), by computing
\begin{equation}
\begin{array}{rcl}
\displaystyle\ln\left(\prod_{k=0}^{p-1}e\left(*\left(n\atop1\right)\right)
(\omega^k z^{1\over p})\right)&=&\\
\displaystyle\sum_{k=0}^{p-1}\ln\left(e\left(*\left(n\atop1\right)\right)
(\omega^k z^{1\over p})\right)&=&\\
 \displaystyle\sum_{k=0}^{p-1}\sum_{j\geq 1}{\left(jn\atop
 j\right)\over j}
\omega^{jk} z^{j\over p}&=&\displaystyle \sum_{k\geq
1}{\left(knp\atop k^p\right)\over
 k}z^k\\
&=&
 \ln(e\left(*\left(np\atop p\right)\right).
\end{array}
\end{equation}

Next (\ref{e1}) follows from the equality
\begin{equation}
e\left(*\left(np\atop p\right)\right)
=
\prod_{k=0}^{p-1}{e\left(\frac1n*\left(n\atop1\right)\right)
(\omega^k z^{1\over p})-1\over \omega^k z^{1\over p}}
\end{equation}
$\Box$

>From \cite{MN} (Theorem 1, p.7), the series $e(*\left(np\atop
p\right))$ is the generating series of the number $q_m$ of lattice
paths from $(0,0)$ to $((n-1)pm,pm)$ that never go above the path
$(\uparrow^{(n-1)p}\rightarrow^{p})^m$
(lattice paths are represented by words over the alphabet
$F=\{\rightarrow,\uparrow\}$,
where $\rightarrow$ means the elementary horizontal path
$(1,0)$ and $\uparrow$ the  elementary vertical path $(0,1)$).

\begin{example}
\rm

The series $$\begin{array}{rl}e(*\left(4\atop 2\right))(z)&=1+6z+53z^2+554z^3+6362z^4+\cdots\\&
=e(*\left(2\atop1\right))(z^{1\over
2})e(*\left(2\atop1\right))(-z^{1\over2})\end{array}$$ is also the
generating series of ordered trees on $2n$ nodes with every
subtree at the root having an even number of edges(see Sloane
\cite{Sloane} {\tt ID number:A066357}). See also \cite{BM} for
another enumeration.

$\begin{array}{rcl}e(*\left(6\atop2\right))(z)&=&1+15z+360z^2+10463z^3+\cdots\\
&=&e(*\left(3\atop 1\right))(z^{1\over 2})e(*\left(3\atop
1\right))(-z^{1\over 2})\end{array}$

$\begin{array}{rcl}e(*\left(6 \atop 3\right))(z)&=&1+20z+662z^2+26780z^3+\cdots\\
&=&e(*\left(2\atop 1\right))(z^{1\over 3})e(*\left(2\atop
1\right))(\exp\{{2i\pi\over 3}\}z^{1\over 3})\times\\&&\times
e(*\left(2\atop 1\right))(\exp\{{4i\pi\over 3}\}z^{1\over
3})\end{array}$

$\begin{array}{rcl}e(*\left(12\atop 3\right))(z)
=&=&1+220z+91498z^2+47961320z^3+\cdots\\
&=&e(*\left(4\atop 1\right))(z^{1\over 3})e(*\left(4\atop
1\right))(\exp\{{2i\pi\over 3}\}z^{1\over 3})\times\\&&\times
e(*\left(4\atop 1\right))(\exp\{{4i\pi\over 3}\}z^{1\over 3})
\end{array}$
\end{example}

One can construct noncommutative analogues of the $*\left(np\atop
p\right)$
by  specializing $S_m$ to the sum of the words coding the
lattice paths from $(0,0)$ to $((n-1)mp,mp)$ that never go above
the path $(\uparrow^{(n-1)p}\rightarrow^p)^m$.

\begin{example}
\rm Let us consider the non-commutative analogue of
$*\left(2\atop 1\right)$. Under this specialization, the first
values of the $S_n$ are

$S_1=\uparrow\rightarrow+\rightarrow\uparrow$

$S_2=\rightarrow^2\uparrow^2+(\rightarrow\uparrow)^2
+\rightarrow\uparrow^2\rightarrow+\uparrow\rightarrow^2\uparrow
+(\uparrow\rightarrow)^2$

$\begin{array}{rcl}S_3
&=&\rightarrow^3\uparrow^3+\rightarrow^2\uparrow\rightarrow\uparrow^2+
\rightarrow^2\uparrow^2\rightarrow\uparrow
+\rightarrow^2\uparrow^3\rightarrow\\
&&+\rightarrow\uparrow\rightarrow^2\uparrow^2
+(\rightarrow\uparrow)^3+(\rightarrow\uparrow)^2\uparrow
\rightarrow+\rightarrow\uparrow^2\rightarrow^2\uparrow\\
&&+\rightarrow\uparrow^2\rightarrow\uparrow\rightarrow
+\uparrow\rightarrow^3\uparrow^2+\uparrow\rightarrow^2
\uparrow\rightarrow\uparrow+\uparrow\rightarrow^2\uparrow^2\rightarrow
+\uparrow\rightarrow\uparrow
\rightarrow^2\uparrow+(\uparrow\rightarrow)^3\end{array}$

\noindent and the first values of the $\Lambda_n$ are

$\Lambda_1=\uparrow\rightarrow+\rightarrow\uparrow$

$\Lambda_2=-\rightarrow^2\uparrow^2$

$\Lambda_3=\rightarrow^2\uparrow\rightarrow\uparrow^2+\rightarrow^3\uparrow^3$
\end{example}
More generally, from
\[\Lambda_m=(-1)^{m+1}\sum_{k=1}^m(-1)^k\sum_{i_1+\cdots+i_k=m}S_{i_1}\cdots S_{i_k}\]
one obtains
\begin{equation}
\Lambda_m=(-1)^{m+1}\sum_{w}w
\end{equation}
where the sum is over  all lattice paths from $(0,0)$ to
$((n-1)mp,mp)$ that never go above
$(\uparrow^{(n-1)p}\rightarrow^{p})^m$ and which avoid the points
$((n-1)kp,kp)$ for $0<k<m$.
Hence,
\begin{equation}
\Phi_m=m\sum_w{w\over c(w)+1}
\end{equation}
where the sum is over the lattice paths from $(0,0)$ to
$((n-1)mp,mp)$ below
$(\uparrow^{(n-1)p}\rightarrow^{p})^m$, and $c(w)$ is the number of
points $((n-1)kp,kp)$ belonging to the path $w$, with $0<k<m$.

\begin{example}\rm
Let us consider the noncommutative analogue of $*\left(4\atop
2\right)$. The first values of the $\Lambda_n$ are

$\Lambda_1=\rightarrow^2\uparrow^2+(\rightarrow\uparrow)^2+\rightarrow\uparrow^2\rightarrow+
\uparrow\rightarrow^2\uparrow+\uparrow\rightarrow\uparrow\rightarrow+\uparrow^2\rightarrow^2$

$\begin{array}{rcl}\Lambda_2&=&-(\rightarrow^4\uparrow^4+\rightarrow^3\uparrow\rightarrow\uparrow^3+
\rightarrow^3\uparrow^2\rightarrow\uparrow^2+\rightarrow^3\uparrow^3\rightarrow\uparrow+
\rightarrow^3\uparrow^4\rightarrow\\ &&+
 \rightarrow^2\uparrow\rightarrow^2\uparrow^3+\rightarrow^2\uparrow(\rightarrow\uparrow)^2\uparrow
  +\rightarrow^2\uparrow\rightarrow\uparrow^2\rightarrow\uparrow+\rightarrow^2\uparrow\rightarrow
  \uparrow^3\rightarrow\\
  &&+\rightarrow\uparrow\rightarrow^3\uparrow^3+\rightarrow\uparrow\rightarrow^2\uparrow\rightarrow\uparrow^2
   +
   \rightarrow\uparrow\rightarrow\uparrow^2\rightarrow\uparrow+\rightarrow\uparrow\rightarrow\uparrow^3\rightarrow\\
  &&
  +\uparrow\rightarrow^4\uparrow^3+\uparrow\rightarrow^3\uparrow\rightarrow\uparrow^2+\uparrow\rightarrow^2
  \uparrow^2\rightarrow\uparrow+\uparrow\rightarrow^2\uparrow^3\rightarrow)\end{array}$
\end{example}
A natural question is whether it is possible to find similar
interpretations for other $*$-binomials.

\subsection{Shuffle analogues of $*$-multinomials in non commutative Witt vectors }

We shall now  construct another noncommutative analogue of
$*$-multinomials. Let $A$ be an alphabet and $w_1,\cdots, w_k$ be
$k$ words of respective lengths $m_1,\ldots,m_k$. Let $w=w_1\cdot
w_2\cdots w_k$ and $m=m_1+\cdots+m_k$. The shuffle $w_1\shuf
\cdots\shuf w_k$ contains ${m\choose m_1,\ldots,m_k}$ terms, and
we can denote it by ${w\choose w_1,\ldots,w_k}$ to emphasize this
point. Then, we introduce the noncommutative Witt vector analogue
\begin{equation}
*{w\choose w_1,\ldots,w_k}=
(Q_1,
Q_2,\cdots,Q_n,\cdots)
\end{equation}
given by the sequence of Witt symmetric functions under the
specialization
\begin{equation}
\Phi_n=w_1^n\shuf w_2^n\shuf\cdots\shuf w_k^n.
\end{equation}
\begin{example}\rm
If $w_i=a^{p_i}$, one has
\begin{eqnarray}\nonumber
*{w\choose w_1,\ldots,w_k}
=(w_1(*\left(p_1+\cdots+p_k\atop
p_1,\dots,p_k\right)a^{p_1+\cdots+p_k},\\\cdots,w_n(*\left((p_1+\cdots+p_k)\atop
p_1,\dots,p_k\right)a^{n(p_1+\cdots+p_k)},\cdots) .\
\end{eqnarray}
Sending $a$ to $1$, one recovers $*$-multinomials.
\end{example}

\section*{Acknowledgements} The authors are grateful to Henri Gaudier
for fruitful discussions and comments.

This project has been partially supported by EC's IHRP Programme, grant
HPRN-CT-2001-00272, ``Algebraic Combinatorics in Europe''.

\end{document}